\documentclass[preprint,3p,times,twocolumn,nopreprintline]{elsarticle} 
\usepackage[utf8]{inputenc}
\usepackage{amsmath}

\usepackage{subcaption}
\usepackage{amsfonts}
\usepackage{adjustbox}
\usepackage{xcolor}
\usepackage{booktabs}
\PassOptionsToPackage{hyphens}{url}\usepackage{hyperref}

\usepackage{natbib}

\usepackage{graphicx}
\usepackage{comment}
\usepackage{appendix}
\usepackage{xspace}
\usepackage{booktabs}
\usepackage{multirow}
\usepackage{makecell}
\usepackage{eurosym}

\PassOptionsToPackage{hyphens}{url}\usepackage{hyperref}
\hypersetup{
    colorlinks
}

\RequirePackage{filecontents}

\usepackage[acronym,toc]{glossaries}

    \newacronym{RREH}{RREH}{Remote Renewable Energy Hub}
    \newacronym{cdo}{CO$_2$}{carbon dioxide}
    \newacronym{PCCC}{PCCC}{Post Combustion Carbon Capture}
    \newacronym{DAC}{DAC}{Direct Air Capture}
    \newacronym{CAPEX}{CAPEX}{Capital Expenditures}
    \newacronym{OPEX}{OPEX}{Operational Expenditures}
    \newacronym{WACC}{WACC}{Weighted Average Cost of Capita}
    \newacronym{gwp}{GWP}{Global Warming Potential}
    \newacronym{ptg}{PtG}{Power to Gas}
    \newacronym{esc}{ESC}{Energy Supply Chain}
    \newacronym{hvdc}{HVDC}{High Voltage DC}

\makeglossaries

\begin{document}
\newcommand{\rreh}{\gls{RREH}\xspace}
\newcommand{\rrehs}{\gls{RREH}s\xspace}
\newcommand{\pccc}{\gls{PCCC}\xspace}
\newcommand{\dac}{\gls{DAC}\xspace}
\newcommand{\cdo}{\gls{cdo}\xspace}
\newcommand{\capex}{\gls{CAPEX}\xspace}
\newcommand{\opex}{\gls{OPEX}\xspace}
\newcommand{\wacc}{\gls{WACC}\xspace}
\newcommand{\gwp}{\gls{gwp}\xspace}
\newcommand{\ptg}{\gls{ptg}\xspace}
\newcommand{\esc}{\gls{esc}\xspace}
\newcommand{\ch}{CH$_4$\xspace}
\newcommand{\hvdc}{\gls{hvdc}\xspace}

\begin{frontmatter}
\title{Synthetic methane for closing the carbon loop: Comparative study of three carbon sources for remote carbon-neutral fuel synthetization}
\author[uliege]{Michaël Fonder}
\author[uliege]{Pierre Counotte}
\author[uliege]{Victor Dachet\corref{cor1}}
\author[tes]{Jehan de Séjournet}
\author[uliege,tpip]{Damien Ernst}
\affiliation[uliege]{organization={University of Liège},
             city={Liège},
             country={Belgium}}
\affiliation[tes]{organization={Tree Energy Solutions},
             city={Zaventem},
             country={Belgium}}
\affiliation[tpip]{organization={Telecom Paris, Institut Polytechnique de Paris},
             city={Paris},
             country={France}}

\cortext[cor1]{Corresponding Author; victor.dachet@uliege.be}
\date{September 2023}
\begin{abstract}
Achieving carbon neutrality is probably one of the most important challenges of the 21st century for our societies. Part of the solution to this challenge is to leverage renewable energies. However, these energy sources are often located far away from places that need the energy, and their availability is intermittent, which makes them challenging to work with. In this paper, we build upon the concept of Remote Renewable Energy Hubs (RREHs), which are hubs located at remote places with abundant renewable energy sources whose purpose is to produce carbon-neutral synthetic fuels. More precisely, we model and study the Energy Supply Chain (ESC) that would be required to provide a constant source of carbon-neutral synthetic methane, also called e-NG (electric Natural Gas) or e-methane (electric methane), in Belgium from an RREH located in Morocco. To be carbon neutral, a synthetic fuel has to be produced from existing carbon dioxide (CO$_2$) that needs to be captured using either Direct Air Capture (DAC) or Post Combustion Carbon Capture (PCCC). In this work, we detail the impact of three different carbon sourcing configurations on the price of the e-methane delivered in Belgium. 
Our results show that sourcing CO$_2$ through a combination of DAC and PCCC is more cost-effective, resulting in a cost of 146€/MWh for e-methane delivered in Belgium, as opposed to relying solely on DAC, which leads to a cost of 158€/MWh. Moreover, these scenarios are compared to a scenario where CO$_2$ is captured in Morocco from a CO$_2$ emitting asset that allow to deliver e-methane for a cost of 136€/MWh.
\end{abstract}

\begin{keyword}
Synthetic Methane \sep Remote Renewable Energy Hub \sep CO2 Sourcing \sep Energy Transition
\end{keyword}

\end{frontmatter}

\section{Introduction}

Global warming is a climate change that is due to large anthropogenic emissions of greenhouse gases, mainly \cdo, in Earth's atmosphere. As its effects are detrimental for our societies and ecosystems in general, it has appeared necessary to reduce, and even cancel, the emission of these gases to minimise these effects. Since most of the \cdo is emitted by burning fossil fuels to generate energy, transitioning our energy production out of these fuels is mandatory.

The main alternatives to fossil fuels are renewable energies such as wind, solar or hydropower, which can be harvested to produce electricity. However, most renewable energy sources are intermittent, and the best sources are located far away from places where this energy is the most needed \cite{Hampp2023Import}. As electrical energy is challenging to transport over long distance and to store in large quantities with current technology, innovations are needed to bridge the gap between energy production and consumption locations.

One innovation to bridge this gap could be the employment of an \hvdc transmission line for the transportation of energy over long distances. Ongoing investigations into projects aiming to link Morocco with the United Kingdom have been documented by \citet{xlinks}. However, this alternative is currently in the early stages of development, and its eventual feasibility remains uncertain.



Another innovation is the concept of \rreh for carbon-neutral fuel synthesis that was first laid out by \citet{Hashimoto1999Global}. The idea underlying \rrehs is to install power-to-X facilities, that convert  electrical energy into chemical energy, at remote locations where renewable energy sources are the most abundant \cite{sadeghi2019energy,geidl2006energy}, and to transport the converted energy to places where it is needed. A power-to-X facility converts the electrical energy into chemical energy by synthesizing energy-dense molecules \cite{LewandowskaBernat2018Opportunities}, which are easier to store and transport back to the energy consumption locations than electricity. The molecules typically considered for this kind of application are dihydrogen, ammonia, or methane, each having its own advantages and drawbacks \cite{Brynolf2018Electrofuels}. 

Among these three molecules, synthetic methane, also called e-methane, holds a particular place. Indeed, it can be used as a simple drop-in replacement within existing energy-demanding installations while decreasing GHG emissions when produced from renewable energy sources, as shown in several case studies \cite{Rixhon2021TheRole,Reiter2015Evaluating,Zhang2017Life}. More generally, existing life-cycle analysis of methanation indicate that the \gwp of e-methane produced with renewable energy sources is smaller than fossil natural gas \cite{Schreiber2020Life,Federici2022Life}. Benefits of e-methane on \gwp are even larger if the methanation is done through catalysts through captured \cdo  rather than from biomass \cite{GoffartDeRoeck2022Comparative,Chauvy2022Environmental}.

To be considered as carbon-neutral, catalyst-based e-methane has to be generated from existing \cdo, which must be actively captured prior to being used for methanation. Capturing \cdo can be done either by filtering it out of ambient air thanks to \dac technologies \cite{McQueen2021AReview,Shayegh2021Future} or by extracting it directly at the source of emission, in industrial fumes, where it is the most concentrated, by \pccc \cite{Zanco2021Postcombustion,Mukherjee2019Review}. Fossil fuel power plants, or cement and steel factories are examples of industries that emit sufficiently large amounts of \cdo to be considered for \pccc. The high concentration of \cdo in industrial fumes leads \pccc to be less expensive than \dac for extracting the same amount of \cdo \cite{Dachet2023Towards-arxiv}. However, capturing \cdo by \pccc needs to be done at the factory that uses the synthesized fuel. The \cdo then needs to be transported where needed, which induces additional costs. This contrasts with \dac that can be performed anywhere on the planet directly where needed, which eliminates the transportation costs. As changing the source of \cdo impacts the price of the final commodity produced by an \rreh, it is necessary to know which \cdo source and capture method combination is the most energy-efficient and/or the most cost-effective for a given \rreh.

Multiple works have studied \ptg \rrehs in different configurations. For example, \citet{Berger2021Remote} studied an \rreh located in Algeria designed to produce e-methane to be delivered in Belgium. \citet{Dachet2023Towards-arxiv} built on this work by adding a hub in Greenland and by studying the impact of carbon pricing on the sizing of the system. \citet{Hampp2023Import} performed an extensive study on the opportunity to build \rrehs at different places on Earth to deliver energy to Germany. As opposed to the works of \citet{Berger2021Remote} and \citet{Dachet2023Towards-arxiv} that focus on e-methane, \citet{Hampp2023Import} consider and compare multiple carriers for energy. However, all these works consider only a single configuration for sourcing the \cdo required for the methanation process in their study. Therefore, it is impossible to know if the \cdo sourcing chosen in each of these works is the most effective one.

This paper aims to compare the impact of changing the \cdo sourcing used to close the carbon loop to the cost of the generated e-methane, and to the sizing of the \rreh, which has not been done before. More specifically, we study the impact of three different configurations of \cdo sources and capture methods on the cost and the sizing of an \rreh located in Morocco and designed to produce e-methane to be delivered in Belgium.

We detail the exact scope and configurations considered in our problem statement, in the next section. We detail our model and methodology in Section 3, and analyse our results in Section 4. Section 5 concludes this paper.

The main contributions of this work are as follows:
\begin{itemize}
    \item We model the whole \esc needed to deliver e-methane, produced from salt water and renewable energy in Morocco, to Belgium;
    \item We consider three different ways of sourcing the \cdo needed for the methanation in the \rreh;
    \item We provide a detailed analysis of the advantages and drawbacks of each of these sources.
\end{itemize}

\section{Problem statement}

In this section, we detail the scope of this work. We consider  an \rreh located in Morocco whose purpose is to produce e-methane for the Belgian market. The \rreh has to produce e-methane from renewable energy sources (for electricity), from salt water and from captured \cdo. The e-methane produced by the hub has to be delivered to Belgium by LNG carriers.  

The choice of Morocco for the location of the hub is first motivated by the quality of its renewable energy sources  \cite{Fasihi2015Economics}.  Second, an \rreh requires a significant land area to be deployed, especially for collecting renewable energies. Morocco is also interesting in this regard, since the southern half of the country is mostly a desert with an extremely low existing land use. 

We consider three different configurations for capturing the \cdo required for the methanation process in the \rreh, namely: (i) \dac on site, in Morocco; (ii) \pccc from Moroccan \cdo emitting asset; and (iii) \pccc from e-methane use in Belgium, with capture losses being compensated by \dac on site, in Morocco.

The purpose of this work is to model the \esc that corresponds to each of these configurations, starting from renewable electricity and salt water in Morocco up to the delivery of e-methane in Belgium, and  to provide an insight on the impact each configuration has on the price of the delivered e-methane and on the sizing of the different elements making the \esc.

The \esc necessary to deliver e-methane to Belgium will be designed to be fully autonomous and auto sufficient. To be comparable to previous works \cite{Berger2021Remote,Dachet2023Towards-arxiv}, the system will be designed to deliver 10TWh of e-methane uniformly over a year in Belgium. The study will be performed with a \wacc set to 7\% , and with extra contingency costs added to all \capex to factor in for feasibility uncertainties~\cite{Roussanaly2021Towards}. More precisely, contingency costs of 10\% and 30\% are to be added to the \capex of mature and unproven technologies respectively.

\section{Material and method}
In this section, we present the method used to carry the study described in the previous section. We first detail the framework used to model and size the desired \esc. We then fix the geographical parameters of our study. Finally, we provide a detailed insight on the setup of the models used  to analyse the three \cdo sourcing configurations considered alongside the parameters used for each part of the model.

\subsection{Modelling framework}
This study is built on the work of \citet{Berger2021Remote} which introduced the use of a Graph-Based Optimization Modeling Language (GBOML) as an optimization framework for multi-energy system models. The GBOML language~\cite{Miftari2022} allows one to model each part of a complex system as a set of nodes interconnected by hyperedges that model the constraints existing between these nodes. In the case of an \rreh each node models a specific module of the hub, and each hyperedge models the flow of a given commodity within the hub.

In the GBOML language, the nodes are, in part, modelled by a set of variables that need to be tuned to minimize the objective function of the model while verifying a set of linear constraints specific to each node. For this, each node has to provide an objective to minimize that is a linear function of its variables. The objective function of the whole model is then defined as the sum of all the objectives of its nodes. In this study, we want to minimize the cost of the e-methane delivered in Belgium, and therefore the cost of the hub. As a result, our objective functions will be proportional to the \capex and \opex of the different modules of the hub.

Since the lifetime of each module is specific, raw \capex are not representative of real costs when studying  the cost of the hub over a limited time horizon. One method to address this concern, which we employ in this study, involves the use of annualized \capex. The annualized \capex $\zeta_m$ of a module $m$ can be computed from the raw \capex, from the life-time $L_m$ of the module and from the \wacc $w$ as follows:
\begin{equation}
\zeta_m = \text{\capex}_m \times \frac{w}{1-(1+w)^{-L_m}}
\end{equation}

Given that the optimization problem of interest has already been comprehensively formalized by \citet{Berger2021Remote}, we shall only remind ourself of the set of assumptions made by prior studies to model and size an \rreh within this framework:
\begin{itemize}
    \item All technologies and their interactions are modelled using linear equations;
    \item The infrastructures and networks needed to transport and process commodities within a single node are not modelled;
    \item Curves modelling the boundaries of the model, such as the renewable energy load factors and the energy demand, are assumed to be known in advance for the whole optimization time horizon;
    \item The sizing of the modules is assumed to be constant over the whole time horizon considered;
    \item Sizing the whole hub with an optimization framework assumes that all planning, investment and operation decisions are made by a single entity.
\end{itemize}

\subsection{Geographical parameters}
\begin{figure*}[t]
    \centering
    \includegraphics[width=0.75\textwidth]{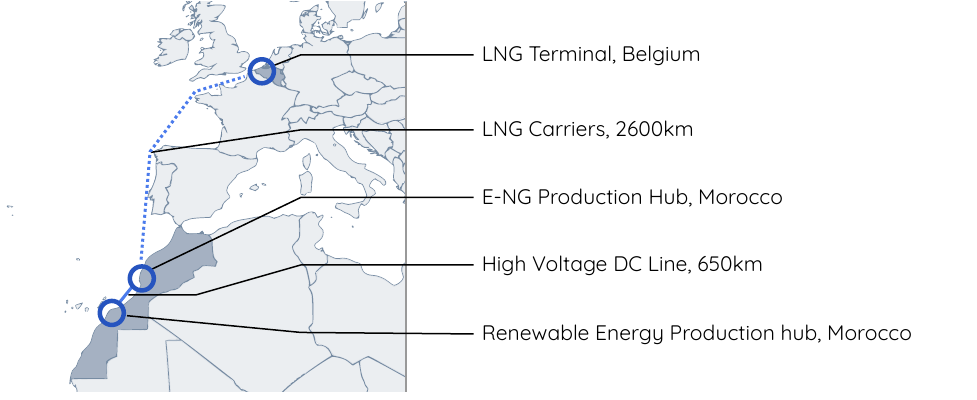}
    \caption{Geographical location of the different parts of the \esc studied in this work.}
    \label{fig:geography}
\end{figure*}
One of the aspects that heavily impact any study of an \rreh is the location of its different components. As mentioned in the problem statement, this study focuses on an \rreh located in Morocco with a e-methane terminal in Belgium. In this part of the document we detail the choice of location of all the parts of the hub and the constraints that result from this choice.

The precise location of the different modules of the hub is mapped in Fig.~\ref{fig:geography}. The delivery terminal is set to be in the Belgian harbour of Zeebrugge as it already features LNG terminals. We build our study on a renewable energy capture hub located on the Atlantic coat in the Western Sahara desert. We consider this location because it has excellent wind and solar energy sources, and because it is mostly empty, which leaves a lot of room to deploy large fields of windmills and solar PV panels. Finally, we set our e-methane production hub to be near the Moroccan city of Safi. By doing so, the e-methane production hub is located near a pool of available workforce, which is required to run the hub, near the coast, which is needed to export the produced e-methane by carriers, and near a large coal-fired power plant~\cite{Nareva2023Safi}, which can provide a valuable source of \cdo for \pccc. It is worth noting that the Safi power plant is sized to deliver 10TWh of electricity over a year~\cite{Nareva2023Safi}, and is therefore a source of \cdo large enough to feed the methanation process of the hub studied in this work. In this study, we consider \pccc applied to the Safi coal-fired power plant. However, our results should be similar to \pccc applied to another \cdo emitting asset, like a steel plant or a cement plant (also present in Safi \cite{HeidelbergCement2023Safi}). In addition, the area of Safi has a lot of available space for large industrial projects, which is a requirement due to the sheer scale of the hub considered. 

With such a setup, the power and e-methane production hubs are separated by 650 km, and need to be connected by a \hvdc line. Since the power production hub is also on the coast, the \hvdc line can be installed offshore to get the shortest connection distance. The e-methane production hub and the LNG terminal in Belgium are separated by 2600 km, and require connection by LNG carriers. We estimate that 100 hours are required to connect Safi to Zeebrugge  by carrier.

\begin{figure*}[t]
    \centering
    \includegraphics[width=0.7\textwidth]{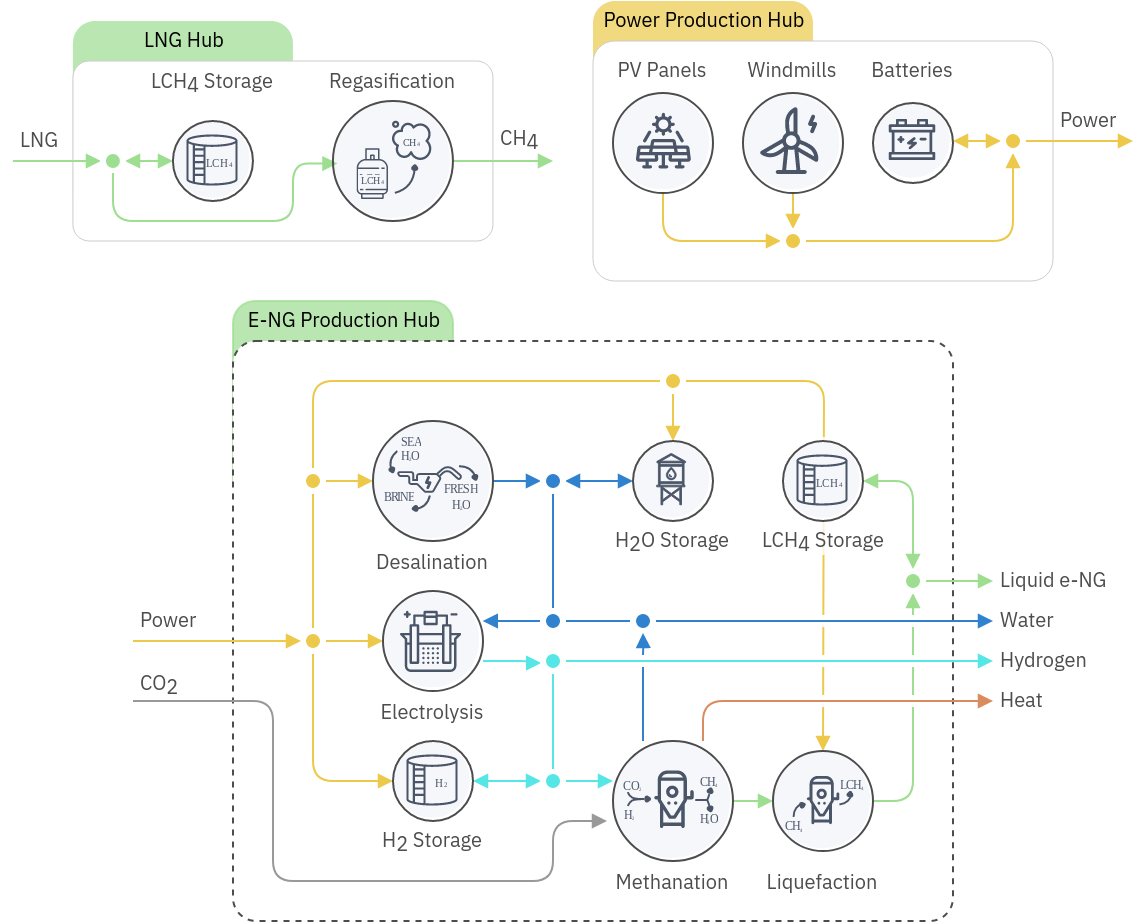}
    \caption{Illustration of three major modules that are part of our \esc and that are common for all \cdo sourcing configurations. The LNG hub features the infrastructures necessary to receive liquefied e-methane from carriers for the end delivery. The module will be installed in Belgium. The power production hub groups all the elements required to provide renewable power coming from energies to other elements of the \esc. Finally, the e-methane production hub takes power and \cdo as inputs and produces liquid e-methane, water, hydrogen, and heat as outputs. These outputs can be consumed by other elements of the \esc, such as \dac units. }
    \label{fig:sub-hubs}
\end{figure*}

\subsection{Model}
Despite being close to the works of \citet{Berger2021Remote} and \citet{Dachet2023Towards-arxiv}, the \rreh model used in this study features some differences. In this section, we first explain the fundamental principles of the hub, and then highlight the key differences with existing works. We would like to emphasize that we open-sourced the code used for the hub to provide all the details required to reproduce this work\footnote{\url{https://gitlab.uliege.be/smart_grids/public/gboml/-/tree/master/examples/synthetic_methane_morocco}}.

\subsubsection{Hub principle}
As detailed in the previous section, the whole \esc necessary to deliver e-methane, synthesized from renewable energy sources in Morocco, in Belgium is made up of the three main hubs illustrated in Fig.~\ref{fig:sub-hubs}. First, the power production hub is made of three nodes; the solar PV panel farm, the windmill farm, and a pack of batteries used to smooth out variations inherent to renewable energy sources. The power production hub is used to deliver the electricity needed to capture \cdo in Morocco, and to synthesize e-methane. The power profiles used to model the renewable energy production over the optimization time horizon were obtained from the \textit{renewable.ninja} website~\cite{Staffell2016Using,Pfenninger2016Longterm,renewablesninja}.

Second, the e-methane production hub, which takes electricity and a source of \cdo as inputs, is responsible for synthesizing e-methane from salt water and \cdo. This hub gets the clean water required for electrolysis from a reverse osmosis water desalination module. This fresh water then passes through an electrolysis module to produce the hydrogen needed by the hub. The core of this hub, the methanation unit, consumes both hydrogen and \cdo to synthesize e-methane with the Sabatier reaction. This node has two byproducts that can be used at other places in the \esc: fresh water and low-grade steam heated to 300°C. Lastly, the synthesized e-methane is liquefied for transportation. In addition, water, hydrogen, and liquefied e-methane storage nodes are added to buffer flow variations of each commodity.

Third, the LNG hub in Belgium is simply made of two nodes: one liquefied e-methane storage and one e-methane regasification node to deliver the e-methane in a gaseous form. 

\begin{figure*}[p]
    \centering
    \includegraphics[width=0.8\textwidth]{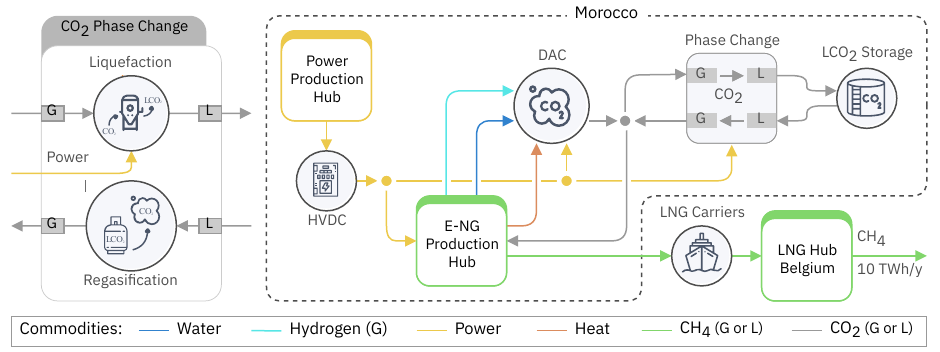}
    \caption{Illustration of the hub design used for our first \cdo sourcing configuration. The \cdo is captured by \dac units at the place where the methanation is performed.}
    \label{fig:sc1}
\end{figure*}

\begin{figure*}[p]
    \centering
    \includegraphics[width=0.6\textwidth]{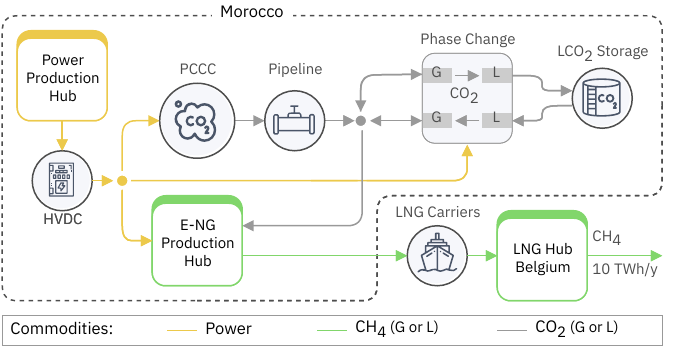}
    \caption{Illustration of the hub design used for our second \cdo sourcing configuration. The \cdo is captured by a \pccc unit placed at a local \cdo emitting asset, in Morocco. The \cdo is transferred to the methanation plant by pipeline.}
    \label{fig:sc2}
\end{figure*}

\begin{figure*}[p]
    \centering
    \includegraphics[width=0.8\textwidth]{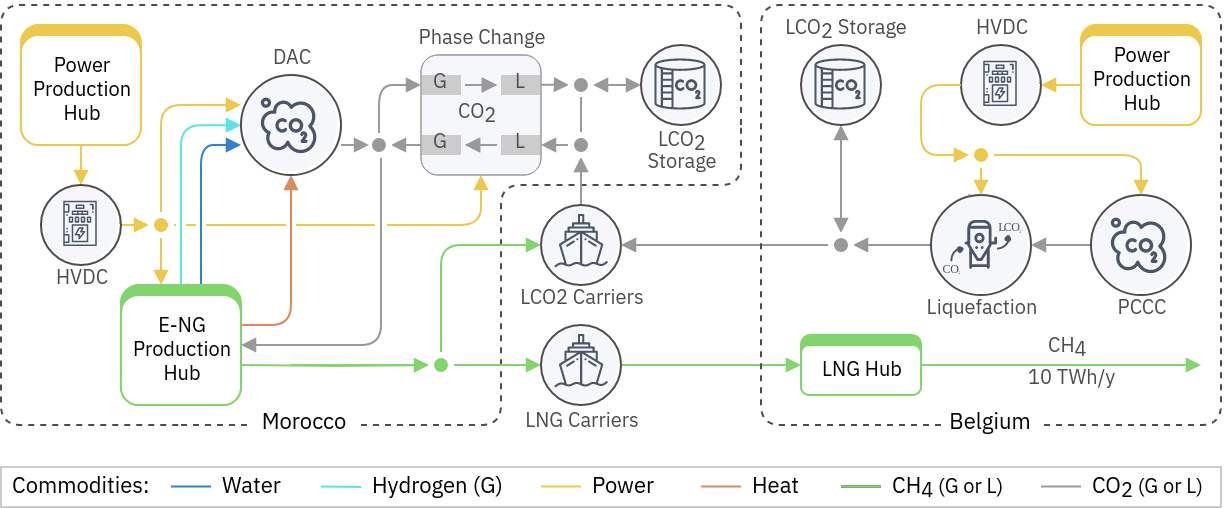}
    \caption{Illustration of the hub design used for our third \cdo sourcing configuration. The \cdo is captured by \pccc at places where e-methane is used in Belgium. This carbon is transferred to the methanation hub by liquefied \cdo carriers. The electricity required for the \pccc units is provided by a small renewable energy power hub. Since \pccc cannot capture all the emitted \cdo, the missing \cdo required for methanation is captured by \dac units on site, in Morocco.}
    \label{fig:sc3}
\end{figure*}

The model of the complete \esc necessary to produce e-methane from \cdo capture by \dac units in Morocco is illustrated in Fig.~\ref{fig:sc1}. The one for \cdo sourced from a \cdo emitting asset \pccc is given in Fig.~\ref{fig:sc2}, while Fig.~\ref{fig:sc3} illustrates the model of the \cdo sourcing configuration implying \pccc in Belgium and \dac in Morocco. In all of these models, \cdo transits as a gas, excepted for storage and carrier transportation where it needs to be liquefied.


\subsubsection{Our adaptations}
Although we use the same nodes with the same parameters as \citet{Berger2021Remote}  and \citet{Dachet2023Towards-arxiv} for most nodes of the model, we need to adapt some parameters of the reference model for our study.  Since  \citet{Berger2021Remote} already detail all the parameters of their model, we only focus on our adaptations to these parameters hereafter.

As required by our problem statement, we add a contingency cost to the \capex of all nodes. This contingency cost, which was not considered in previous works, is equal to 10\% of the \capex for well-established technologies, and 30\% of the \capex for less mature technologies. The technologies that we consider to be less mature regarding the scale of the \rreh are the following: large-scale battery packs, \hvdc lines, electrolysis units, methanation units, and \cdo capture technologies. 

Modifications are also carried out to the methanation node when compared to the original node proposed by \citet{Berger2021Remote}. First, we align our \capex with the  $300\text{k€/}(\text{MWH \ch/h})$ proposed by  \citet{Gorre2019Production}. Second, we add an output for the residual heat produced by the Sabatier reaction as it can also be valued. Indeed,  \citet{Coppitters2023Energy} showed that this heat can be recycled in the hub to lower the  heat energy required by solid sorbent \dac units. According to their work, the Sabatier reaction produces 2.1MWh of usable heat energy per ton of synthesized \ch.

Using a heat source of 300°C  to regenerate a sorbent is only possible for solid sorbents. Indeed, solid sorbents only require to be heated to 100°C to regenerate~\cite{Coppitters2023Energy}, whereas liquid sorbents need to be heated to 900°C~\cite{Berger2021Remote}. Therefore, as opposed to preceding works~\cite{Berger2021Remote,Dachet2023Towards-arxiv}, we use solid sorbent \dac units in our model, and rely on the data of the \citet{DanishEnergyAgency} for their specifications. In addition to using recirculated heat, the system has the ability to burn hydrogen to provide heat for the sorbent regeneration.

A last minor change of in our model when compared to previous works is the cost of liquefaction, storage and regasification of \cdo. The source used in previous works~\cite{Mitsubishi2004Ship} dates from 2004, which is significantly older than most of our sources. To uniformize cost estimates, we adjusted the estimates of \citet{Mitsubishi2004Ship} by taking inflation between 2004 and 2021 into account.

All the techno-economic parameters for each modeled technology are listed in the tables available in the Appendices.

\begin{table*}[t]
\centering
\begin{tabular}{l*{4}{r}}
\toprule
\makecell[l]{ECS Elements} &  \makecell[r]{\ch Production\\Free \cdo \\(€/MWh \ch)} & \makecell[r]{ \cdo Sourcing\\DAC MA\\(€/MWh \ch)} & \makecell[r]{ \cdo Sourcing\\PCCC MA\\(€/MWh \ch)} &  \makecell[r]{\cdo Sourcing \\PCCC BE + DAC MA\\(€/MWh \ch)} \\ 
\midrule 
Solar PV Field & 10.43 & 1.23 & 0.37 & 0.11 \\ 
Windmill farm & 43.00 & 4.02 & 1.9 & 0.20 \\ 
Batteries & 4.06 & 1.49 & -0.05 & 0.33 \\ 
HVDC & 11.86 & 1.12 & 0.50 & 0.06 \\ 
\cdo capt. + transport. & 0.00 & 21.61 & 8.39 & 20.09 \\ 
Desalination plant & 0.11 & 0.66 & 0.00 & 0.16 \\ 
Water storage & 0.05 & 0.00 & 0.00 & 0.00 \\ 
Electrolysis plant & 35.42 & 2.89 & 0.29 & 0.05 \\ 
H2 storage & 4.55 & 0.44 & 0.17 & 0.02 \\ 
Methanation plant & 7.75 & 0.00 & 0.00 & 0.01 \\ 
\ch liquefaction & 5.09 & 0.00 & 0.00 & 0.01 \\ 
\ch storage (Morocco) & 0.23 & 0.00 & 0.00 & 0.00 \\ 
\ch carriers & 0.67 & 0.00 & 0.00 & 0.00 \\ 
\ch storage (Belgium) & 0.23 & 0.00 & 0.00 & 0.00 \\ 
\ch regasification & 1.02 & 0.00 & 0.00 & 0.00 \\ 
\midrule  
Total cost & 124.47 & 33.46 & 11.57 & 21.04 \\ 
\midrule  
\ch cost  & 124.47 & 157.93 & 136.04 & 145.51 \\ 
\bottomrule
\end{tabular}
\caption{Breakdown of the cost for delivering 1MWh of \ch in Belgium depending on the source of \cdo used for the methanation process. The cost given for the "\cdo capt. + transport." line encompasses all the infrastructures required to capture and transport the \cdo back to the e-methane production hub. This line is detailed in Table~\ref{tab:co2-costs} for each \cdo sourcing configuration individually. The other costs induced by a \cdo sourcing, which listed in the table, detail the costs of feeding the appropriate commodities (electricity, water,...)  to a given \cdo sourcing configuration with the hub.}
\label{tab:costs}
\end{table*}

\section{Results}
With the study and the model being defined, we can analyse the results of the optimization process. In this section, we first analyse the three \cdo sourcing configurations from a cost perspective. We then analyze them from an energy perspective. Finally, we discuss the results obtained in this work. 

\subsection{Costs analysis}
The cost for delivering \ch to Belgium for our three \cdo sourcing configurations are detailed in Tables~\ref{tab:costs} and~\ref{tab:co2-costs}. The price for the commodities produced by the hub is given in Table~\ref{tab:commodity-costs}. Producing and delivering 1MWh of \ch costs 124.47€ without factoring in costs related to \cdo sourcing. The largest part of this cost is due to the renewable energy capture (53.43€/MWh \ch), electricity transportation (11.86€/MWh \ch), and water electrolysis (35.42€/MWh \ch).

Capturing and transporting \cdo to the e-methane production hub creates additional costs that range from 11.57€/MWh \ch for the \pccc in Morocco to 33.46€/MWh \ch for \dac on site, which is the most expensive way to source \cdo to the hub. Therefore, delivering 1MWh of e-methane in Belgium costs a total of 136.04€ when sourcing\cdo from \pccc in Morocco, 145.51€ when sourcing \cdo from \pccc in Belgium, and 157.93€ when sourcing \cdo from \dac on site. For comparison, regular fossil natural gas has been exchanged at prices well above these costs for several months during the 2022 energy crisis in Europe (up to 342€/MWh in August 2022).

\cdo coming from \dac on site seems to also be more expensive than \pccc \cdo imported from Belgium. As a result,  the methanation uses as much \cdo from Belgium as possible in the configuration where \cdo can be sourced either from \dac on site or from \pccc in Belgium. This leads the fraction of \cdo coming from Belgium used for methanation to be equal to the efficiency of the \pccc process (90\% in our model).

In order to be exhaustive, we tested an alternate model where \ch and \cdo are transported through offshore pipelines for the configuration where part of the \cdo comes from \pccc in Belgium. In this case, the price of e-methane delivered in Belgium rises to 151.24€/MWh, which is more expensive than transportation by carriers.

\begin{table*}[t]
    \begin{minipage}{0.45\textwidth}
        \begin{subtable}[h!]{\linewidth}
            \centering
            \begin{tabular}{l*{1}{r}}
            \toprule
                \makecell[l]{DAC MA\\\cdo sourcing elements} & \makecell[r]{Cost\\(€/Mwh)} \\ 
                \midrule 
                DAC (Morocco) & 21.61 \\ 
                \cdo storage (Morocco) & 0.0 \\ 
                \midrule
                Total cost & 21.61 \\ 
                \bottomrule
                \end{tabular}
            \caption{\dac in Morocco\vspace{0.5cm}}
            \label{tab:sub-dac}
        \end{subtable}\\
        
        \begin{subtable}[h!]{\linewidth}
            \centering
            \begin{tabular}{l*{1}{r}}
            \toprule
                \makecell[l]{PCCC MA\\\cdo sourcing elements} & \makecell[r]{Cost\\(€/Mwh)} \\
                \midrule 
                PCCC (Morocco) & 8.08 \\ 
                \cdo storage (Morocco) & 0.3 \\ 
                \cdo pipe & 0.01 \\ 
                \midrule 
                Total cost & 8.39 \\ 
                \bottomrule
            \end{tabular}
            \caption{\pccc from a \cdo emitting asset in Morocco}
            \label{tab:sub-pccc}
        \end{subtable}
    \end{minipage}
    \hfill
    \begin{minipage}{0.45\textwidth}
        \begin{subtable}[h]{\textwidth}
            \centering
            \begin{tabular}{l*{1}{r}}
            \toprule
                \makecell[l]{PCCC BE + DAC MA\\\cdo sourcing elements} & \makecell[r]{Cost\\(€/Mwh)} \\
                \midrule 
                Solar PV field (Belgium) & 1.12 \\ 
                Windmill farm (Belgium) & 2.36 \\ 
                Batteries (Belgium) & 1.16 \\ 
                HVDC (Belgium) & 0.24 \\ 
                PCCC (Belgium) & 7.07 \\ 
                \cdo liquefaction (Belgium) & 0.14 \\ 
                \cdo storage (Belgium) & 0.44 \\ 
                \cdo carriers & 1.77 \\ 
                DAC (Morocco) & 5.28 \\ 
                \cdo storage (Morocco) & 0.45 \\ 
                \cdo regasification (Morocco) & 0.06 \\ 
                \midrule 
                Total cost & 20.09 \\ 
                \bottomrule
            \end{tabular}
            \caption{\pccc in Belgium completed by \dac in Morocco}
            \label{tab:sub-pcccdac}
        \end{subtable}
    \end{minipage}
    \caption{Breakdown of the cost for sourcing \cdo to the methanation process for our three different sourcing configurations. The electricity required for the PCCC in Belgium is provided by a small renewable energy hub installed in Belgium.}
    \label{tab:co2-costs}
\end{table*}

\begin{table}[ht]
\centering
\begin{tabular}{l*{1}{r}}
\toprule
Commodity & Cost [€] \\ 
\midrule 
Hydrogen [t] & 3210.69  \\ 
Water [t] & 1.17 \\
\cdo - \dac MA [t] & 256.28 \\ 
\cdo - \pccc MA [t] & 72.28 \\ 
\cdo - \pccc BE + \dac MA [t] & 136.62 \\ 
Wind power - MA [MWh] & 30.19 \\ 
Solar power - MA [MWh] & 22.89 \\ 
Wind power - BE [MWh] & 65.47 \\ 
Solar power - BE [MWh] & 58.51 \\ 
\bottomrule
\end{tabular}
\caption{Cost of the commodities produced by our \esc model. Due to limitations, the price given for \cdo sources using \dac correspond to the upper bound price that happens when all the heating energy comes from hydrogen. The production costs for power, water, and hydrogen are independent of the \cdo sourcing configuration.}
\label{tab:commodity-costs}
\end{table}

\subsection{Energy analysis}
\begin{table*}[ht]
\centering
\begin{tabular}{l*{4}{r}}
\toprule
\makecell[l]{\cdo Sourcing\\Configuration} & \makecell[r]{\esc - Global\\Efficiency [\%]} & \makecell[r]{Wind - Fraction of\\Curtailment [\%]}  & \makecell[r]{PV - Fraction of\\Curtailment [\%]} & \makecell[r]{DAC - Part of\\Recycled Heat[\%]} \\ 
\midrule 
DAC MA & 48.44 & 25.17 & 0.0 & 65.41 \\ 
PCCC MA & 51.10 & 25.17 & 0.0 & N/A \\ 
DAC MA + PCCC BE & 51.26 & 25.14 & 0.0 & 99.99\\ 
\bottomrule
\end{tabular}
\caption{Insight of some key figures related to energy within our model. The first column gives the ratio between the energy delivered in Belgium in the form of \ch and the captured renewable energy. The second column gives the fraction of available wind energy that was not used by the model. The third column gives the fraction of available solar energy that was not used by the model. The last column gives the fraction of heat energy required by \dac that comes from recirculated heat rather than from hydrogen.}
\label{tab:energy-stats}
\end{table*}
Table~\ref{tab:energy-stats} gives some key statistics related to energy within our model. This table shows that the efficiency of the hub, that is the ratio between the energy contained in the e-methane delivered in Belgium and the energy captured by windmills and solar panels, lies around 50\% with small variations depending on the \cdo sourcing configuration. The highest efficiency, 51.26\%, is obtained by the \cdo sourcing that combines \pccc in Belgium and \dac in Morocco and beats the efficiency of simple \pccc in Morocco. This observation, which may be surprising at first glance, can be explained by the use of heat recirculation. As shown in the last column of the table, a large part of the heat required for \dac can be provided by the Sabatier reaction. In the case of joint \pccc and \dac, the Sabatier reaction can even provide all the heat required for \dac, which lowers the energy input required for carbon capture when compared to \pccc alone. However, it is worth noting that the \cdo sourcing involving \pccc in Morocco has residual heat that is not used in our model, but that could be valued for other applications, which would virtually increase the efficiency of the hub.

The Table~\ref{tab:energy-stats} also shows the fraction of available renewable energy that is not used by the hub. For windmills, the energy curtailment rises to 25\%, while all the available solar energy is used. This observation is identical for all \cdo sourcing configurations. We believe that the discrepancy between the curtailments of solar and wind power can be explained by a difference in regularity of the energy sources. Indeed, the available solar power is pretty regular from a day to the other whereas wind power fluctuates more.

In addition to these insights, we computed the load factor for different parts of the system. Windmills and PV panels have a load factor of 41.4\% and 25.5\% respectively. Electrolysis units have a load factor of about 81\%, while desalination and methanation units are used at full capacity all the time by design. These load factors are similar for all \cdo sourcing configurations.

Finally, Table~\ref{tab:sizing} shows the capacity that needs to be installed for delivering an average power of 1MW of e-methane to Belgium for various modules of the hub. Since the model is fully linear, these capacities scale linearly with the desired power of e-methane to be delivered. It is interesting to note that the biggest consumers of fresh water are by far\dac units according to numbers given in this table.
\begin{table*}[t]
\centering
\begin{tabular}{l*{4}{r}}
\toprule
\makecell[l]{ECS Elements} &  \makecell[r]{\ch Production\\Free \cdo \\ {[}/MW \ch]} & \makecell[r]{ \cdo Sourcing\\DAC MA\\{[}/MW \ch]} & \makecell[r]{ \cdo Sourcing\\PCCC MA\\ {[}/MW \ch]} &  \makecell[r]{\cdo Sourcing \\PCCC BE + DAC MA\\{[}/MW \ch]} \\ 
 \midrule 
PV Capacity [MW] &     2.1389 &     0.2514 &     0.0769 &     0.0226 \\ 
 Windmill capacity [MW] &     3.4660 &     0.3238 &     0.1528 &     0.0162 \\ 
 Battery capacity [MWh] &     1.1266 &     0.4182 &    -0.0188 &     0.0909 \\ 
 Battery throughput [MW] &     0.1469 &     0.0469 &     0.0007 &     0.0110 \\ 
 HVDC line capacity [MW] &     2.3105 &     0.2178 &     0.0983 &     0.0110 \\ 
 Desalination units [MW] &     0.1508 &     0.9468 &     0.0000 &     0.2257 \\ 
 Electrolysis units [MW] &     2.1717 &     0.1774 &     0.0179 &     0.0029 \\ 
 \bottomrule
 \end{tabular}
\caption{Capacity  of several modules of the hub required to deliver an average power of 1MW of e-methane in Belgium.}
\label{tab:sizing}
\end{table*}

\subsection{Results discussion}
All the results presented in this section are directly linked to our hypothesis and to our model. Our observations and analysis may be slightly different when considering some additional factors. In this section, we discuss how some factors that could impact the results presented in this work.

First, our model does not factor in for the economies of scale that could come with the deployment of a large \rreh. Indeed, most of the current cost estimated are projections based on existing prototypes or units that are of small scale or still at experimental stage. Therefore, it seems reasonable to assume  that the overall cost of modules would fall if manufactured in large quantities.

Second, the efficiency of \pccc can vary depending on the industry on which \cdo capture is performed~ \cite{Reiter2015Evaluating}. As little data are available on the type and capacity of industries that may use the delivered e-methane in Belgium, we had to make an educated assumption on the efficiency of \pccc. Real-world deployments may lead to efficiencies different to the assumptions made in this work, and could change our observations. Indeed, if the efficiency of \pccc drops to steeply, costs associated to \pccc may rise to a point where the balance between \pccc in Belgium and \dac in Morocco has a different optimum.

Third, we modelled the hub to be self-sufficient. As shown in the previous section, this implies to synthesize hydrogen to provide heat to \dac units when the heat resulting from the Sabatier reaction is not sufficient, which increases the cost of the hub. However, the low-grade heat required to regenerate solid sorbents in \dac units is a common industrial byproduct that is often considered as a waste due to a lack usecase, but that can find a use within the hub. As a result, importing heat from  a neighbouring industry could make the price of the configuration relying on \dac units more competitive price-wise when compared to other configurations of \cdo sourcing.

Fourth, the ability of \dac units to provide a \cdo source that does not need a dedicated supply chain is a practical advantage over \pccc that does not appear when simply looking at numbers. A simpler supply chain for the production of e-methane is also probably a more robust one as fewer elements come into play. This could lead to a higher reliability of the whole supply chain, which has an economic value that is not taken into account in our model.

Lastly, we need to discuss the carbon neutrality of the synthetic methane in the scenario implying \pccc on an existing \cdo emitting asset. While the other scenarios are clearly carbon neutral due to the complete recapture of \cdo emitted by the use of e-methane, performing \pccc on an existing \cdo emitting asset can only be considered to be carbon neutral if this asset was installed prior to or is operated independently from the \rreh. As assets have a limited lifetime, 35 years in the case of a coal-fired power plant \cite{coal_lifetime}, this scenario could be seen as an intermediate step towards one of the two others scenarios since it offers the opportunity to get an operating \rreh in a short term and under a cost-effective budget. Furthermore, in all scenarios, emissions of \cdo in Belgium are covered by \cdo withdrawal or avoided emissions. Therefore, the carbon neutrality of e-methane should be recognized by the exchange of emission permits between Belgium and Morocco, which in turn will incentivize the use of e-methane in Belgium.


\section{Conclusion}

In this paper, we consider the sizing and the cost of an \esc designed to synthesize e-methane from renewable energy sources in an \rreh in Morocco and to deliver it in Belgium. Synthesizing e-methane requires a source of \cdo. In this work, we considered  three different configurations for sourcing the required \cdo to the \rreh, and studied their impact on the sizing and the cost of the \rreh, and by extension on the cost of the e-methane delivered in Belgium. The three  different configurations considered for capturing the \cdo are (i) \dac on site, in Morocco; (ii) \pccc from Moroccan \cdo emitting asset; and (iii) \pccc from e-methane use in Belgium, with capture losses being compensated by \dac on site, in Morocco.

We modelled and optimized the \esc corresponding to the three configurations using the GBOML framework. Results show that \dac is more expensive than \pccc generally speaking. This is in part due to the heat required to regenerate the sorbent of \dac units. Indeed, the \cdo sourcing relying only on \dac units leads to a price of  157.93€/MWh of e-methane delivered in Belgium, while the configuration relying on \pccc from Moroccan \cdo emitting asset lead to of price of 136.04€/MWh.  The hybrid configuration, which implies \pccc in Belgium and \dac in Morocco, leads to an intermediate price of 145.51€/MWh of e-methane delivered in Belgium. The price for \cdo induced by this configuration is still lower than the one of \cdo captured by \dac on site despite the additional costs induced by the transportation of \cdo between Belgium and Morocco. These costs, which are estimates based on conservative projections on technology costs by the year 2030, can be expected to decrease with every technological breakthroughs linked to elements of the \rreh.



\noindent\textbf{Acknowledgements.~}
 Victor Dachet was supported by the Walloon Region (Service Public de Wallonie Recherche, Belgium) under grant n°2010235– ARIAC by \textit{digitalwallonia4.ai}.

\printglossaries

\bibliographystyle{elsarticle-num-names.bst}
\bibliography{abbreviation,abbreviation-short,refs}

\section*{Appendices}\label{sec:appendices}

The tables (\ref{tab:stor_econ_stock_params}, \ref{tab:stor_econ_flow_params}, \ref{tab:stor_tech_params}, \ref{tab:conv_tech_params}, \ref{tab:conv_econ_params}) have been adapted from \citet{Berger2021Remote} to encompass all the technologies used in this study. Economical costs of storage components are presented in two separate tables—one considering the flow and the other focusing on the stock of commodities. For readers interested in a comprehensive mathematical formulation of the constraints and objective functions related to these parameters, we recommend referring to \citet{Berger2021Remote}.

\begin{table*}
	
	\centering
	
	\begin{tabular}{lccccc}
		\hline
		& CAPEX & FOM & VOM & Lifetime & Contingency \\ 
		\hline
		Battery Storage & \small{142.0} & \small{0.0} & \small{0.0018} & \small{10.0} & 0.3 \\
		 \citep{batterycosts} & \small{M\euro/GWh} & \small{M\euro/GWh-yr} & \small{M\euro/GWh} & \small{yr}  & - \\
		Compressed H$_2$ Storage & \small{45.0} & \small{2.25} & \small{0.0} & \small{30.0} & 0.1 \\
		 \citep{batterycosts} & \small{M\euro/kt} & \small{M\euro/kt-yr} & \small{M\euro/kt} & \small{yr} & - \\
		Liquefied CH$_4$ Storage & \small{2.641} & \small{0.05282} & \small{0.0} & \small{30.0} & 0.1 \\
		 \citep{lngtankcosts} & \small{M\euro/kt}& \small{M\euro/kt-yr}& \small{M\euro/kt} & \small{yr} &  - \\
        Liquefied CO$_2$ Storage & \small{2.3} & \small{0.0675} & \small{0.0} & \small{30.0} & 0.1 \\
		\citep{co2liquefaction} & \small{M\euro/(kt/h)} & \small{M\euro/(kt/h)} & \small{M\euro/kt} & \small{yr} & - \\
		H$_2$O Storage & \small{0.065} & \small{0.0013} & \small{0.0} & \small{30.0} & 0.1 \\
		 \citep{Caldera2016} & \small{M\euro/kt} & \small{M\euro/kt-yr} & \small{M\euro/kt} & \small{yr} & - \\
		\hline
	\end{tabular}
 \caption{Economic parameters used to model storage nodes (stock component, 2030 estimates).}
	\label{tab:stor_econ_stock_params}
\end{table*}

\begin{table*}
	\centering

	\begin{tabular}{lccccc}
		\hline
		& CAPEX & FOM & VOM & Lifetime & Contingency \\ 
		\hline
		Battery Storage & \small{160.0} & \small{0.5} & \small{0.0} & \small{10.0}  & 0.3 \\
		 \citep{batterycosts} & \small{M\euro/GW} & \small{M\euro/GW-yr} & \small{M\euro/GWh} & \small{yr} & - \\
		Compressed H$_2$ Storage & \small{45.0} & \small{2.25} & \small{0.0} & \small{30.0} & 0.1 \\
		 \citep{batterycosts} & \small{M\euro/kt} & \small{M\euro/kt-yr} & \small{M\euro/kt} & \small{yr} & - \\
		Liquefied CH$_4$ Storage & \small{2.641} & \small{0.05282} & \small{0.0} & \small{30.0} & 0.1 \\
		 \citep{lngtankcosts} & \small{M\euro/kt}& \small{M\euro/kt-yr}& \small{M\euro/kt} & \small{yr} &  - \\
        Liquefied CO$_2$ Storage & \small{ 0.0} & \small{ 0.0} & \small{ 0.0} & \small{30} & 0.1 \\
		\citep{co2liquefaction} & \small{M\euro/(kt/h)} & \small{M\euro/(kt/h)} & \small{M\euro/kt} & \small{yr} & - \\
		H$_2$O Storage & \small{1.55923} & \small{0.0312} & \small{0.0} & \small{30.0} & 0.1 \\
		 \citep{Caldera2016} & \small{M\euro/(kt/h)} & \small{M\euro/(kt/h)} & \small{M\euro/kt} & \small{yr} & - \\
		\hline
	\end{tabular}
 \caption{Economic parameters used to model storage nodes (flow component, 2030 estimates).}
	\label{tab:stor_econ_flow_params}
\end{table*}

\begin{table*}
	
	\centering
	\begin{tabular}{lcccccc}
		\hline
		&  \small{$\eta^S$} & \small{$\eta^+$} & \small{$\eta^-$} & \small{$\sigma$} & \small{$\rho$} & \small{$\phi$} \\ \hline
		Battery Storage & \small{0.00004} & \small{0.959} & \small{0.959} & \small{0.0} & \small{1.0} & \\
		\citep{batterycosts} & \small{-} & \small{-} & \small{-} & \small{-} & \small{-} &\\
		Compressed H$_2$ Storage & \small{1.0} & \small{1.0} & \small{1.0} & \small{0.05} & \small{1.0} & \small{1.3}\\
		\citep{batterycosts} & & & & & & \small{GWh$_{el}$/kt$_{H_2}$}\\
		Liquefied CO$_2$ Storage & \small{1.0} & \small{1.0} & \small{1.0} & \small{0.0} & \small{1.0} & \small{0.105}\\
		\citep{co2liquefaction} & & & & & & \small{GWh$_{el}$/kt$_{CO_2}$}\\
		Liquefied CH$_4$ Storage & \small{1.0} & \small{1.0} & \small{1.0} & \small{0.0} & \small{1.0} & \\
		 & \small{-} & \small{-} & \small{-} &\small{-} & \small{-} & \\
		H$_2$O Storage & \small{1.0} & \small{1.0} & \small{1.0} & \small{0.0} & \small{1.0} & \small{0.00036}\\
		\citep{Caldera2016} & & & & & &\small{GWh$_{el}$/kt$_{H_2O}$}\\
		\hline
	\end{tabular}
 \caption{Technical parameters used to model storage nodes. The $\eta^S$ corresponds to the self discharge rate, $\eta^+$ the charge efficiency, $\eta^-$ the discharge efficiency, $\sigma$ the minimum inventory level, $\rho$ the maximum discharge-to-charge ratio and $\phi$ the conversion factor.}
	\label{tab:stor_tech_params}
\end{table*}

\begin{table*}[h]
	\begin{tabular}{lcccccc}
	\hline
	& $\phi_1$ & $\phi_2$ & $\phi_3$ & $\phi_4$ & $\mu$ & $\Delta_{+,-}$ \\
	\hline
	HVDC Interconnection & \small{0.9499} & & & & & \\
	\citep{Xiang2016,hvdclineefficiency} & \small{-} & & & & & \\
	Electrolysis & \small{50.6} & \small{9.0} & \small{8.0} & & \small{0.05} & \small{1.0}\\
	\citep{Goetz2016} & \small{GWh$_{el}$/kt$_{H_2}$} & \small{kt$_{H_2O}$/kt$_{H_2}$} & \small{kt$_{O_2}$/kt$_{H_2}$} & & \small{-} & \small{-/h} \\
	Methanation & \small{0.5} & \small{2.75} & \small{2.25} & \small{0,1345} &  \small{1.0} & \small{0.0} \\
	\citep{Goetz2016,Roensch2016} & \small{kt$_{H_2}$/kt$_{CH_4}$} & \small{kt$_{CO_2}$/kt$_{CH_4}$} & \small{kt$_{H_2O}$/kt$_{CH_4}$} & \small{GWh$_{heat}$/GWh$_{CH_4}$} & \small{-} & \small{-/h}\\
	Desalination & \small{0.004} & & & & \small{1.0} & 0.0\\
	\citep{desalinationconsumption} & \small{GWh$_{el}$/kt$_{H_2O}$} & & & &  \small{-} & \small{-/h}\\
	Direct Air Capture & \small{ 0.15 } & \small{5.0} & \small{1.46} & \small{0.2} & \small{1.0} & \small{0.0} \\
	   \citep{DanishEnergyAgency} &\small{GWh$_{el}$/kt$_{CO_2}$}&\small{kt$_{H_2O}$/kt$_{CO_2}$} &\small{GWh$_{heat}$/kt$_{CO_2}$}& \small{GWh$_{heat}$/kt$_{CO_2}$} & \small{-} & \small{-/h} \\
    Post Combustion Carbon Capture & \small{ 0.4125} &  & \small{ } & \small{} & \small{0.0} & \small{1.0} \\
	   \citep{DanishEnergyAgency} &\small{GWh$_{el}$/kt$_{CO_2}$} &  &  & & \small{-} & \small{-/h} \\
	CH$_4$ Liquefaction & \small{0.616} & & & &  \small{0.0} & 1.0\\
	\citep{Pospisil2019} & \small{GWh$_{el}$/kt$_{LCH_4}$} & & & &  \small{-}  & \small{-/h}\\
	LCH$_4$ Carriers & \small{0.994} & & &  & &\\
	\citep{lngcarrierlifetime} & \small{-} & & & & & \\
	LCH$_4$ Regasification & \small{0.98} & & & & &\\
	\citep{Pospisil2019} & \small{-} & & & & &\\
    CO$_2$ Liquefaction & \small{0.014} & .99 & & &  \small{0.0} & 1.0\\
	\citep{DanishEnergyAgency} & \small{GWh$_{el}$/kt$_{LCO_2}$} & & & &  \small{-}  & \small{-/h}\\
	LCO$_2$ Carriers & \small{0.99} & \small{0,0000625} & &  & &\\
	\citep{DanishEnergyAgency} & \small{-} & \small{(GWh/h)/kt$_{CO_2}$} & & & & \\
	LCO$_2$ Regasification & \small{0.98} & & & & &\\
	\citep{DanishEnergyAgency} & \small{-} & & & & &\\
    Pipe CO$_2$ & \small{0.00002} & & & & &\\
	\citep{DanishEnergyAgency} & \small{GWh/kt$_{CO_2}$/h} & & & & &\\
	\hline
	\end{tabular}
 \caption{Technical parameters used to model conversion nodes. The $\phi_i$ represents the conversion factors, $\mu$ the minimum operating level , $\Delta_{+,-}$ the maximum ramp-up/ramp-down rate.}
	\label{tab:conv_tech_params}
\end{table*}

\begin{table*}

	\resizebox{\textwidth}{!}{%
	\begin{tabular}{lccccc}
	\hline
	& CAPEX & FOM & VOM & Lifetime & Contingency \\
	\hline
	Solar Photovoltaic Panels & \small{380.0} & \small{7.25} & \small{0.0} & \small{25.0}  & 0.1 \\
	\citep{solarpvcosts} & \small{M\euro/GW$_{el}$} & \small{M\euro/GW$_{el}$-yr} & \small{M\euro/GWh$_{el}$} & \small{yr} & - \\
	Wind Turbines & \small{1040.0} & \small{12.6} & \small{0.00135} & \small{30.0} & 0.1\\
	\citep{windcosts} & \small{M\euro/GW$_{el}$} & \small{M\euro/GW$_{el}$-yr} & \small{M\euro/GWh$_{el}$} & \small{yr} & - \\
	HVDC Interconnection & \small{480.0} & \small{7.1} & \small{0.0} & \small{40.0} & 0.3 \\
	\citep{hvdccost,hvdccostsEIA} & \small{M\euro/GW$_{el}$} & \small{M\euro/GW$_{el}$-yr} & \small{M\euro/GWh$_{el}$} & \small{yr} & - \\
	Electrolysis & \small{600.0} & \small{30.0} & \small{0.0} & \small{15.0} & 0.3 \\
	\citep{electrolysiscosts} & \small{M\euro/GW$_{el}$} & \small{M\euro/GW$_{el}$-yr} & \small{M\euro/GWh$_{el}$} & \small{yr} & - \\
	Methanation & \small{300.0} & \small{29.4} & \small{0.0} & \small{20.0}  & 0.3 \\
	\citep{Goetz2016, methanationcosts} & \small{M\euro/GW$_{CH_4}$ (HHV)} & \small{M\euro/GW$_{CH_4}$-yr (HHV)} & \small{M\euro/GWh$_{CH_4}$ (HHV)} & \small{yr} & - \\
	Desalination & \small{28.08} & \small{0.0} & \small{0.000315} & \small{20.0}  & 0.1 \\
	\citep{desalinationcosts} & \small{M\euro/(kt$_{H_2O}$/h)} & \small{M\euro/(kt$_{H_2O}$/h)-yr} & \small{M\euro/kt$_{H_2O}$} & \small{yr} & - \\
	Direct Air Capture & \small{6000.0} & \small{300.0} & \small{0.0} & \small{20.0}  & 0.3 \\
	 \citep{DanishEnergyAgency}   &\small{M\euro/(kt$_{CO_2}$/h)} &\small{M\euro/(kt$_{CO_2}$/h)-yr} &\small{M\euro/kt$_{CO_2}$} & \small{yr} & - \\
    Post Combustion Carbon Capture & \small{3150.0} & \small{0.0} & \small{0.0} & \small{20.0}  & 0.3 \\
	\citep{DanishEnergyAgency}   &\small{M\euro/(kt$_{CO_2}$/h)} &\small{M\euro/(kt$_{CO_2}$/h)-yr} &\small{M\euro/kt$_{CO_2}$} & \small{yr} & - \\
	CH$_4$ Liquefaction & \small{5913.0} & \small{147.825} & \small{0.0} & \small{30.0} & 0.1 \\
	\citep{liquefactioncosts} & \small{M\euro/(kt$_{LCH_4}$/h)} & \small{M\euro/(kt$_{LCH_4}$/h)-yr} & \small{M\euro/kt$_{LCH_4}$} & \small{yr} & - \\
	LCH$_4$ Carriers & \small{2.537} & \small{0.12685} & \small{0.0} & \small{30.0} & 0.1  \\
	\citep{lngcarriercapex} & \small{M\euro/kt$_{LCH_4}$} & \small{M\euro/kt$_{LCH_4}$-yr} & \small{M\euro/kt$_{LCH_4}$} & \small{yr} & - \\
	LCH$_4$ Regasification & \small{1248.3} & \small{24.97} & \small{0.0} & \small{30.0} & 0.1 \\
	\citep{regasificationcosts} & \small{M\euro/(kt$_{CH_4}$/h)} &\small{M\euro/(kt$_{CH_4}$/h)-yr} &\small{M\euro/kt$_{CH_4}$} & \small{yr} & - \\
    CO$_2$ Liquefaction & \small{55.8} & \small{2.79} & \small{0.0} & \small{30.0} & 0.1 \\
	\citep{DanishEnergyAgency} & \small{M\euro/(kt$_{LCO_2}$/h)} & \small{M\euro/(kt$_{LCO_2}$/h)-yr} & \small{M\euro/kt$_{LCO_2}$} & \small{yr} & - \\
	LCO$_2$ Carriers & \small{5} & \small{0.0} & \small{0.0} & \small{40.0} & 0.1 \\
	\citep{DanishEnergyAgency} & \small{M\euro/kt$_{LCO_2}$} & \small{M\euro/kt$_{LCO_2}$-yr} & \small{M\euro/kt$_{LCO_2}$} & \small{yr} & - \\
	LCO$_2$ Regasification & \small{25.1} & \small{1.25} & \small{0.0} & \small{30.0} & 0.1 \\
	\citep{Mitsubishi2004Ship} & \small{M\euro/(kt$_{CO_2}$/h)} &\small{M\euro/(kt$_{CO_2}$/h)-yr} &\small{M\euro/kt$_{CO_2}$} & \small{yr} & - \\
    Pipe CO$_2$  & \small{2.3} & \small{20.0} & \small{0.0} & \small{40.0} & 0.1 \\
	\citep{DanishEnergyAgency} & \small{M\euro/kt$_{CO_2}$/km} & \small{M\euro/kt$_{CO_2}$-yr} & \small{M\euro/kt$_{CO_2}$} & \small{yr} & - \\
	\hline
	\end{tabular}}
 \caption{Economic parameters used to model conversion nodes (2030 estimates).}
	\label{tab:conv_econ_params}
\end{table*}

\end{document}